\documentclass[a4paper,11pt]{article}
\usepackage{latexsym}
\usepackage{amssymb}
\usepackage{amsfonts}
\usepackage{amsmath}
\usepackage{amsthm}
\usepackage{MnSymbol}
\usepackage{calligra}
\usepackage[boxruled]{algorithm2e}
\usepackage[T1]{fontenc}
\usepackage{mathrsfs}
\usepackage{indentfirst}
\usepackage{graphicx}
\usepackage{epsfig}
\usepackage[enableskew]{youngtab}
\usepackage{color}
\usepackage{fullpage}
\usepackage{tikz,ifthen}
\newcommand{\pione}{\pi_1}

\newtheorem{prop}{Proposition}[section]
\newtheorem{teor}{Theorem}[section]
\newtheorem{lemma}{Lemma}[section]
\newtheorem{cor}{Corollary}[section]

\newcommand{\ninN}{n\in \mathbf{N}}

\newcommand{\cvd}{\hfill $\blacksquare$\bigskip}

\newcommand{\Sch}{${\mathscr S}\hspace{-4pt}${\calligra ch}\; \,}
\newcounter{indice}

\newcommand{\shape}[1]{
\setcounter{indice}{0}; \foreach \i in {#1} {
\addtocounter{indice}{1};
\foreach \x in {1,...,\i} { \draw
(\x-1,-\theindice+1) rectangle (\x,-\theindice); \draw
(\x,-\theindice+1) -- (\x-1,-\theindice);} } }

\newcommand{\secondshape}[1]{
\setcounter{indice}{0}; \foreach \i in {#1} {
\addtocounter{indice}{1};
\foreach \x in {1,...,\i} { \draw
(\x+5,-\theindice+1) rectangle (\x+6,-\theindice); \draw
(\x+6,-\theindice+1) -- (\x+5,-\theindice);} } }

\newcommand{\Sshape}[1]{

}

\date{}

\author{Luca Ferrari\thanks{ Dipartimento di Matematica e
Informatica ``U. Dini'',, Universit\`a degli Studi di Firenze,
Viale Morgagni 65, 50134 Firenze, Italy,
\texttt{luca.ferrari@unifi.it}. Partially supported by INdAM -
GNCS 2015 project ``Problemi di consistenza, unicit\'a e
ricostruzione per grafi ed ipergrafi'' and by MIUR PRIN 2010-2011 grant ``Automi e
Linguaggi Formali: Aspetti Matematici e Applicativi'', code
H41J12000190001.}}

\title{Schr\"oder partitions, Schr\"oder tableaux and weak poset patterns}

\begin{document}

\maketitle

\begin{abstract}
We introduce the notions of \emph{Schr\"oder shape} and of \emph{Schr\"oder tableau},
which provide some kind of analogs of the classical notions of Young shape and Young tableau.
We investigate some properties of the partial order given by containment
of Schr\"oder shapes. Then we propose an algorithm which is the natural analog of the
well known RS correspondence for Young tableaux, and we characterize those permutations
whose insertion tableaux have some special shapes. The last part of the article relates the notion
of Schr\"oder tableau with those of interval order and of weak containment (and strong avoidance)
of posets. We end our paper with several suggestions for possible further work.
\end{abstract}

\section{Introduction}

Given a positive integer $n$, a \emph{partition} of $n$ is a
finite sequence of positive integers $\lambda =(\lambda_1
,\lambda_2 ,\ldots ,\lambda_r )$ such that $\lambda_1 \geq
\lambda_2 \geq \cdots \geq \lambda_r$ and $n=\lambda_1 +\lambda_2
+\cdots +\lambda_r$. When $\lambda$ is a partition of $n$ we also
write $\lambda \vdash n$. A graphical way of representing
partitions is given by Young shapes. The \emph{Young shape} of the
above partition $\lambda \vdash n$ consists of $r$ left-justified
rows having $\lambda_1 ,\ldots ,\lambda_r$ boxes (also called
cells) stacked in decreasing order of length. The set of all Young
shapes can be endowed with a poset structure by containment (of
top-left justified shapes). Such a poset turns out to be in fact a
lattice, called the \emph{Young lattice}. A \emph{standard Young
tableau} with $n$ cells is a Young shape whose cells are filled in
with positive integers from $1$ to $n$ in such a way that entries
in each row and each column are (strictly) increasing.

Young tableaux are among the most investigated combinatorial
objects. The widespread interest in Young tableaux is certainly
due both to their intrinsic combinatorial beauty (which is
witnessed by several surprising facts concerning, for instance,
their enumeration, such as the hook length formula \cite{FRT} and the RSK
algorithm \cite{R,Sc,Kn}) and to their usefulness in several algebraic contexts,
typically in the representation theory of groups and related
matters (such as Schur functions and the Littlewood-Richardson
rule \cite{LR}).

Apart from their classical definition, there are several
alternative ways to introduce Young tableaux. In the present paper
we are interested in the possibility of defining standard Young
tableaux in terms of a certain lattice structure on Dyck paths.
The main advantage of this point of view lies in the possibility
of giving an analogous definition in a modified setting, in which
Dyck paths are replaced by some other class of lattice paths. Here
we will try to see what happens if we replace Dyck paths with
Schr\"oder paths, just scratching the surface of a theory that, in
our opinion, deserves to be better studied.

\bigskip

Given a Cartesian coordinate system, a \emph{Dyck path} is a
lattice path starting from the origin, ending on the $x$-axis,
never falling below the $x$-axis and using only two kinds of
steps, $u(p)=(1,1)$ and $d(own)=(1,-1)$. A Dyck path can be
encoded by a word $w$ on the alphabet $\{ u,d\}$ such that in
every prefix of $w$ the number of $u$'s is greater than or equal to
the number of $d$'s and the total number of $u$ and $d$ in $w$ is
the same (the resulting language is called \emph{Dyck language}
and its words \emph{Dyck words}). The \emph{length} of a Dyck path
is the length of the associated Dyck word (which is necessarily an
even number).

\begin{figure}[h]
\begin{center}
 \setlength{\unitlength}{2mm}
 \begin{picture}(100,15)
% \put(0,0){
% \setlength{\unitlength}{1.5mm}
% \begin{picture}(15,18)
%  \put(0,0){\vector(1,0){17}}
%  \put(0,0){\vector(0,1){19}}
%  \multiput(0,0)(1,3){7}{\circle*{0.3}}
%  \multiput(2,1)(1,3){6}{\circle*{0.3}}
%  \multiput(4,2)(1,3){5}{\circle*{0.3}}
%  \multiput(5,0)(1,3){5}{\circle*{0.3}}
%  \multiput(7,1)(1,3){4}{\circle*{0.3}}
%  \multiput(9,2)(1,3){3}{\circle*{0.3}}
%  \multiput(10,0)(1,3){3}{\circle*{0.3}}
%  \multiput(12,1)(1,3){2}{\circle*{0.3}}
%  \put(14,2){\circle*{0.3}}
%  \put(15,0){\circle*{0.3}}
%  \put(0,0){\line(1,3){6}}
%  \put(2,1){\line(1,3){5}}
%  \put(4,2){\line(1,3){4}}
%  \put(5,0){\line(1,3){4}}
%  \put(7,1){\line(1,3){3}}
%  \put(9,2){\line(1,3){2}}
%  \put(10,0){\line(1,3){2}}
%  \put(12,1){\line(1,3){1}}
%  \put(2,1){\line(-1,2){1}}
%  \put(5,0){\line(-1,2){3}}
%  \put(7,1){\line(-1,2){4}}
%  \put(10,0){\line(-1,2){6}}
%  \put(12,1){\line(-1,2){7}}
%  \put(15,0){\line(-1,2){9}}
% \end{picture}}
% \qbezier(12,7)(13,8.7)(14,9)
% \put(14,9){\vector(4,1){0.3}}
 \put(14.5,3){
 \setlength{\unitlength}{2.6mm}
 \begin{picture}(15,10)(0,-3.5)
  \put(0,0){\vector(1,0){15}}
  \put(0,-3.5){\vector(0,1){10}}
  \multiput(0,0)(1,1){7}{\circle*{0.2}}
  \multiput(3,1)(1,1){3}{\circle*{0.2}}
  \put(7,5){\circle*{0.2}}
  \put(6,2){\circle*{0.2}}
  \put(8,4){\circle*{0.2}}
  \put(7,1){\circle*{0.2}}
  \put(9,3){\circle*{0.2}}
  \multiput(8,0)(1,1){3}{\circle*{0.2}}
  \put(11,1){\circle*{0.2}}
  \put(12,0){\circle*{0.2}}
  \thicklines
  \put(0,0){\color{red}\line(1,1){2}}
  \put(2,2){\line(1,1){4}}
  \put(3,1){\color{red}\line(1,1){2}}
  \put(8,0){\color{red}\line(1,1){2}}
  \thinlines
  \put(3,1){\color{red}\line(-1,1){1.5}}
  \thicklines
  \put(5,3){\color{red}\line(1,-1){3}}
  \put(6,6){\line(1,-1){4}}
  \put(10,2){\color{red}\line(1,-1){2}}
  \color{white}
  \put(1.9,2.1){\line(-1,1){1.5}}
  \color{black}
  \thinlines
  \put(5,3){\line(1,1){2}}
  \put(6,2){\line(1,1){2}}
  \put(7,1){\line(1,1){2}}
  \put(3,3){\line(1,-1){1.5}}
  \put(4,4){\line(1,-1){1.5}}
  \put(5,5){\line(1,-1){4}}
  \color{white}
  \thicklines
  \put(4,2){\line(1,-1){1.5}}
  \color{black}
  \put(4,2){\circle*{0.2}}
  \end{picture}}
 \qbezier(33,9)(34,9.9)(35,10)
 \put(35,10){\vector(4,1){0.3}}
 \put(36,8){
 \setlength{\unitlength}{4mm}
 \begin{picture}(9,6)
  %\put(0,0){\dashbox{0.2}(9,6){}}
%  \multiput(0,0)(0,1){7}{\circle*{0.2}}
%  \multiput(1,1)(0,1){6}{\circle*{0.2}}
%  \multiput(2,2)(0,1){5}{\circle*{0.2}}
%  \multiput(3,2)(0,1){5}{\circle*{0.2}}
%  \multiput(4,3)(0,1){4}{\circle*{0.2}}
%  \multiput(5,4)(0,1){3}{\circle*{0.2}}
%  \multiput(6,4)(0,1){3}{\circle*{0.2}}
%  \multiput(7,5)(0,1){2}{\circle*{0.2}}
%  \put(8,6){\circle*{0.2}}
%  \put(9,6){\circle*{0.2}}
  \put(0,0){\line(0,1){4}}
  \put(1,0){\line(0,4){4}}
  \put(2,2){\line(0,1){2}}
  \put(3,2){\line(0,1){2}}
  \put(4,2){\line(0,1){2}}
  \put(0,0){\line(1,0){1}}
  \put(0,1){\line(1,0){1}}
  \put(0,2){\line(1,0){4}}
  \put(0,3){\line(1,0){4}}
  \put(0,4){\line(1,0){4}}
  \end{picture}}
 %\put(41,6.7){\makebox(0,0){$\downarrow$}}
% \put(36,0){
% \setlength{\unitlength}{2.4mm}
% \begin{picture}(7,5)
%  \put(0,0){\line(0,1){5}}
%  \put(1,0){\line(0,1){5}}
%  \put(2,1){\line(0,1){4}}
%  \put(3,1){\line(0,1){4}}
%  \put(4,2){\line(0,1){3}}
%  \put(5,3){\line(0,1){2}}
%  \put(6,3){\line(0,1){2}}
%  \put(7,4){\line(0,1){1}}
%  \put(0,0){\line(1,0){1}}
%  \put(0,1){\line(1,0){3}}
%  \put(0,2){\line(1,0){4}}
%  \put(0,3){\line(1,0){6}}
%  \put(0,4){\line(1,0){7}}
%  \put(0,5){\line(1,0){7}}
%  \put(0.2,1.2){\dashbox{0.2}(2.6,1.6){}}
%  \put(0.2,3.2){\dashbox{0.2}(2.6,1.6){}}
%  \put(3.2,3.2){\dashbox{0.2}(2.6,1.6){}}
% \end{picture}}
 \end{picture}
\end{center}
 \vspace{-1cm}
 \caption{A Dyck path (red) and the associated Young shape.}
 \label{Dyckshape}
\end{figure}
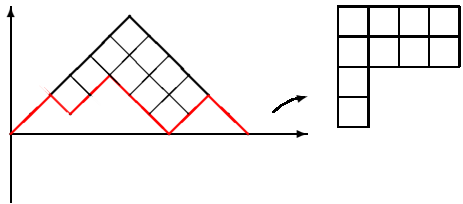

Consider the set $\mathbf{D}_n$ of all Dyck paths of length
$2n$; it can be endowed with a very natural poset structure, by
declaring $P\leq Q$ whenever $P$ lies weakly below $Q$ in the
usual two-dimensional drawing of Dyck paths (for any $P,Q\in
\mathbf{D}_n$). This partial order actually induces a distributive
lattice structure on $\mathbf{D}_n$, to be denoted $\mathcal{D}_n$
and called \emph{Dyck lattice of order $n$}. This can be shown
both in a direct way, using the combinatorics of lattice paths
(see \cite{FP}), and as a consequence of the fact that
$\mathcal{D}_n$ is order-isomorphic to (the dual of) the Young
lattice of the staircase partition $(n-1,n-2,\ldots ,2,1)$ (that
is the principal down-set generated by such a staircase partition
in the Young lattice). Referring to the latter approach, any $P\in
\mathbf{D}_n$ uniquely determines a Young shape, which can be
obtained by taking the region included between $P$ and the maximum
path of $\mathcal{D}_n$, then slicing it into square cells using
diagonal lines of slope $1$ and $-1$ passing through all points
having integer coordinates, and finally rotating the sheet of
paper by $45^{\circ}$ anticlockwise (see Figure \ref{Dyckshape}).

It is well known that there is a bijection between standard Young
tableaux of a given shape and saturated chains in the Young
lattice starting from the empty shape and ending with that shape.
Translating this fact on Dyck lattices, we can thus state that
standard Young tableaux of a given shape are in bijection with
saturated chains (inside a Dyck lattice of suitable order)
starting from the Dyck path associated with that shape and ending
with the maximum of the lattice. This suggests us to try to find an
analog of this fact in which Dyck paths are replaced by other
types of paths. As already mentioned, the case treated in the
present paper is that of Schr\"oder paths.

\bigskip

In section \ref{poset} we introduce the notion of Schr\"oder shape
and study some properties of the poset of Schr\"oder shapes
(in some sense analogous to those of the Young lattice).
In section \ref{algorithm} we introduce the notion of Schr\"oder tableau
and we define an algorithm which, given a permutation,
produces a pair of Schr\"oder tableaux having the same Schr\"oder shape;
this is made in analogy with the classical RS algorithm. In particular,
we will address the problem of determining which permutations are mapped
into the same Schr\"oder insertion tableau, and we solve it for a few special shapes.
Section \ref{alternative} offers an alternative description of the notion of
Schr\"oder tableau in terms of two seemingly unrelated concepts: one is well known
(interval orders) whereas the other one (weak pattern poset, and strong poset avoidance)
is much less studied; we then give an overview of a possible combinatorial approach
to the study of weak poset containment and strong poset avoidance, and provide a link
between these notions and Schr\"oder tableaux. Finally, we devote Section \ref{conclusions}
to the presentation of some directions of further research.

\bigskip

An extended abstract of the present work has appeared in the proceedings of the conference
IWOCA 2015 \cite{Fe}.

\section{The poset of Schr\"oder partitions}\label{poset}

\begin{figure}
\begin{center}
\begin{tikzpicture}[scale=0.5]
\shape{4,3,3,1}

\draw (4,0) -- (5,0) -- (4,-1) -- cycle;

\draw (1,-3) -- (2,-3) -- (1,-4) -- cycle;

\draw (0,-4) -- (1,-4) -- (0,-5) -- cycle;
\end{tikzpicture}
\end{center}
\caption{A Schr\"oder shape of order 25.\label{shape}}
\end{figure}
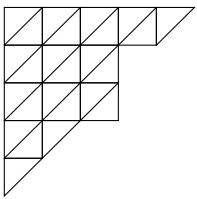

A \emph{Schr\"oder shape} is a set of triangular cells in the
plane obtained from a Young shape by drawing the NE-SW diagonal of
each of its (square) cells, and possibly adding below the first column and at the end of some
rows one more triangular cell, provided that, in a group of rows having equal length,
only the first (topmost) one can have an added triangle.
The number of cells of a Schr\"oder shape is called the \emph{order} of that
shape. An example of a Schr\"oder shape is illustrated in Figure
\ref{shape}.

A Schr\"oder shape has triangular cells of two distinct types,
which will be referred to as \emph{lower triangular cells} and
\emph{upper triangular cells}. In particular, rows having an odd
number of cells necessarily terminate with an upper triangular
cell. A Schr\"oder shape determines a unique integer partition,
whose parts are the number of cells in the rows of the shape. For
instance, the partition associated with the shape in Figure
\ref{shape} is $(9,6,6,3,1)$. As a consequence of the definition
of a Schr\"oder shape, it is clear that not every partition can be
represented using a Schr\"oder shape. More precisely, we have the
following result, whose proof is completely trivial and so it is
left to the reader.

\begin{prop} An integer partition can be represented with a
Schr\"oder shape if and only if its odd parts are \emph{simple}
(i.e. have multiplicity 1).
\end{prop}

Those integer partitions which can be represented with a suitable
Schr\"oder shape will be called \emph{Schr\"oder partitions}. The
set of all Schr\"oder partitions will be denoted $\mathbf{Sch}$,
and the set of Schr\"oder partitions of order $n$ with
$\mathbf{Sch}_n$. From now on we will frequently refer to
Schr\"oder shapes and to Schr\"oder partitions interchangeably,
when no confusion is likely to arise.

From the enumerative point of view, the number of Schr\"oder
partitions is known, and is recorded in \cite{Sl} as sequence
A006950. In particular, the generating function of Schr\"oder
partitions is given by
$$
\prod_{k>0}\frac{1+x^{2k-1}}{1-x^{2k}}.
$$

There are several combinatorial interpretations for the resulting
sequence, however an appropriate reference for the present one (in
terms of Schr\"oder partitions) appears to be \cite{D}. In that
paper the author proves a far more general result, concerning
partitions such that the multiplicity of each odd part is in a
prescribed set and the multiplicity of each even part is
unrestricted.

It is interesting to notice that this sequence is also relevant
from an algebraic point of view. Indeed it coincides with the
sequence of numbers of nilpotent conjugacy classes in the Lie
algebras $o(n)$ of skew-symmetric $n\times n$ matrices. This
suggests that Schr\"oder partitions have a role in representation
theory that certainly deserves to be better investigated.

\bigskip
%
%Here we propose a refined enumerative result, namely we describe a
%simple recurrence for the number of Schr\"oder partitions of $n$
%into $k$ parts.
%
%\begin{prop} Denote with $s_{n,k}$ the number of Schr\"oder partitions of $n$
%into $k$ parts and with $s'_{n,k}$ the number of Schr\"oder
%partitions of $n$ into $k$ parts having smallest part different
%from 1. Then, for all $n\geq k\geq 1$:
%\begin{itemize}
%
%  \item[(i)] $s_{n,k}=s'_{n,k}+s'_{n-1,k-1}$;
%
%  \item[(ii)] $s'_{n,k}=s'_{n-2,k-1}+s'_{n-2k-1,k-1}+s'_{n-2k,k}$.
%
%\end{itemize}
%\end{prop}
%
%\emph{Proof.}\quad We immediately observe that the set of
%Schr\"oder partitions of $n$ into $k$ parts whose smallest part is
%equal to 1 is in bijection with the set of Schr\"oder partitions
%of $n$ into $k-1$ parts whose smallest part is different from 1.
%This gives at once the formula in (i).
%
%Concerning (ii), given a a Schr\"oder partition $\lambda$ of $n$ into $k$ parts
%with no part equal to $1$, we distinguish two cases.
%If $\lambda$ has at least one part equal to 2,
%then removing it leaves us with a Schr\"oder partition of $n-2$ into $k-1$ parts,
%still having no part equal to 1. Otherwise, removing the first two columns of $\lambda$
%returns a Schr\"oder partition of $n-2k$ into $k$ parts, possibly having one part equal to 1.
%From here, using (i), we immediately obtain (ii). \cvd

Though the formalism of Schr\"oder shapes seems not to add relevant
information on the enumerative combinatorics of Schr\"oder partitions,
it suggests at least an interesting family of maps on integer partitions,
which turns out to define a family of involutions if suitably restricted.
Consider the family of maps $(c_n )_{\ninN}$
defined on the set of all integer partitions as follows:
given a partition $\lambda$ and a positive integer $n$,
$c_n (\lambda )$ is the integer partition $\mu =(\mu_1 ,\mu_2 ,\ldots ,\mu_k )$
(of the same size as $\lambda$) whose $i$-th part $\mu_i$ is given by
the sum of the $n$ columns of (the Young shape of) $\lambda$
from the $((i-1)n+1)$-th one to the $(in)$-th one. So, for instance,
$c_3 ((7,6,6,6,4,3,3,1))=(22,13,1)$. Since each of the above maps preserves
the size of a partition, it is clearly an endofunction when restricted to
the set of all integer partitions of size $n$. Notice that $c_1$ is the well-known conjugation map
(which exchanges rows with columns in a Young shape).
In spite of the fact that $c_1$ is an involution (on the set of all partitions),
it is easy to see that all the other $c_n$'s are not involutions.
However, it is possible to characterize the set of those partitions
for which $c_n ^2$ acts as the identity map. The next proposition, which was incorrectly formulated
in \cite{Fe}, is now stated and proved in a correct way.

\begin{prop}
Given $\ninN$ and an integer partition $\lambda$
(whose $i$-th part will be denoted $\lambda_i$, as usual), we have that
$c_n ^2 (\lambda )=\lambda$ if and only if for all $k\geq 0$, there exists at most one index $i$
such that $kn<\lambda_i <(k+1)n$.
\end{prop}

\emph{Proof.}\quad For any given $\lambda$, suppose that there exists at least one part of $c_n (\lambda )$
which is not multiple of $n$, and let $\mu'$ be one of them. More precisely, let $k$ be the unique nonnegative integer
such that $kn<\mu'<(k+1)n$. This means that $\lambda$ has a set of $n$ consecutive columns whose sum is equal to $\mu '$.
Since $\mu '\not \equiv 0\; (\textnormal{mod } n)$, this implies that such $n$ columns are not all equal.
In particular, the rightmost of them must have $<k+1$ cells.
Now, since in a Young shape columns are in decreasing order of length,
the sum of the successive $n$ columns of $\lambda$ is $\leq kn$,
hence $\mu '$ is the only part of $c_n (\lambda )$ strictly greater than $kn$.
Using a similar argument, we observe that, in the set of columns of $\lambda$ that sum up to $\mu'$,
the first (leftmost) of them must have $\geq k+1$ cells, and so the sum of the previous $n$ columns of $\lambda$ is
$\geq (k+1)n$; as a consequence, $\mu'$ is the only part of $c_n (\lambda )$ strictly smaller than $(k+1)n$.
We have thus proved that the condition in the above statement holds for every partition in the image of $c_n$.
This is enough to conclude that, if $c_n ^2 (\lambda )=\lambda$, then necessarily
the same condition holds for $\lambda$ (which lies indeed in the image of $c_n$).

\begin{figure}[!h]
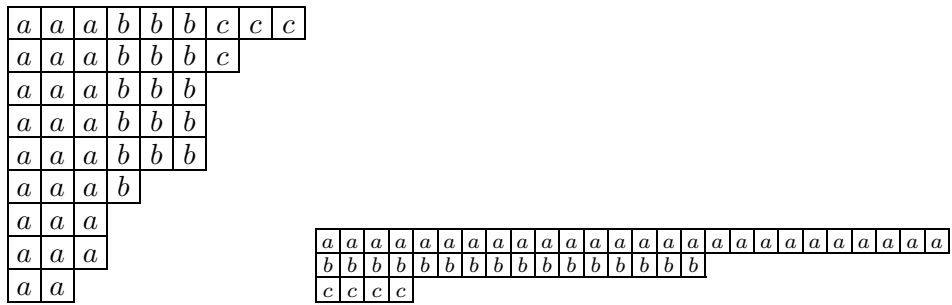

\begin{center}
\young(a:a:a:b:b:b:c:c:c,a:a:a:b:b:b:c,a:a:a:b:b:b,a:a:a:b:b:b,a:a:a:b:b:b,a:a:a:b,a:a:a,a:a:a,a:a)
{\scriptsize \young(a:a:a:a:a:a:a:a:a:a:a:a:a:a:a:a:a:a:a:a:a:a:a:a:a:a,b:b:b:b:b:b:b:b:b:b:b:b:b:b:b:b,c:c:c:c)}
\end{center}
\caption{From $\lambda$ to $c_n (\lambda )$.\label{proposizione}}
\end{figure}

Conversely, split each row of $\lambda$ into clusters containing $n$ consecutive cells,
except at most the last cluster which contains at most $n$ cells.
Denoting with $\mu_i$ the $i$-th part of $c_n (\lambda )$, we have that $\mu_i$ is obtained by
%summing a set of $n$ columns of $\lambda$, namely from column $(i-1)n+1$ to column $in$. The hypothesis on $\lambda$
%implies that the contribution of each row of $\lambda$ to $\mu_i$ is equal to $n$ except at most one row,
%whose contribution is strictly less than $n$ (clearly excluding those rows whose contribution is 0).
%In other words, $\mu_i$ is obtained by
taking the $i$-th cluster from each row, and the hypothesis implies that,
among the rows whose contribution is nonzero, there is at most one row whose contribution is strictly less than $n$.
The construction of $c_n (\lambda )$ from $\lambda$ is illustrated in Figure \ref{proposizione} for the partition
$\lambda =(9,7,6,6,6,4,3,3,2)$ and $n=3$: cells with the same label have to be grouped together,
and the resulting partition $c_3 (\lambda )=(26,16,4)$ is depicted on the right.

Now, similarly as above, in order to construct $c_n ^2 (\lambda )$, we have to split each row of $c_n (\lambda )$ into clusters.
Notice that, as a consequence of our hypothesis,
if the $i$-th row of $c_n (\lambda )$ has a (necessarily unique) cluster containing strictly less than $n$ cells,
then this is precisely the unique cluster with less than $n$ cells among all the $i$-th clusters of all rows of $\lambda$.
Therefore, constructing $c_n ^2 (\lambda )$ from $c_n (\lambda )$ we recover exactly the starting partition $\lambda$, as desired.
%As a consequence, when forming the partition $c_n ^2 (\lambda )$ from $c_n (\lambda )$,
%the $i$-th part is obtained by taking
%Conversely, observe that we can represent every partition $\lambda$
%by means of a Young-like shape, which is obtained from the usual Young shape of $\lambda$
%by simply grouping together the cells of each row $n$ by $n$.
%In this way we obtain a shape (call it $\tilde{\lambda}$)
%in which each cell is a horizontal rectangle made of $n$ cells of the original Young shape,
%except at most the last cell of each row,
%which is a horizontal rectangle having \emph{at most} $n$ cells.
%Now observe that $c_n (\lambda )$ can be obtained by exchanging the rows and the columns
%of $\tilde{\lambda}$ and then breaking the horizontal rectangles of the resulting shape
%into $n$ square cells. This construction is illustrated below for the partition
%$\lambda =(9,7,6,6,6,4,3,3,2)$ and $n=3$: cells with the same label have to be grouped together,
%and the resulting partition $c_3 (\lambda )=(26,16,4)$ is depicted on the right.
%
%$$
%\young(a:a:a:b:b:b:c:c:c,a:a:a:b:b:b:c,a:a:a:b:b:b,a:a:a:b:b:b,a:a:a:b:b:b,a:a:a:b,a:a:a,a:a:a,a:a)
%{\scriptsize \young(a:a:a:a:a:a:a:a:a:a:a:a:a:a:a:a:a:a:a:a:a:a:a:a:a:a,b:b:b:b:b:b:b:b:b:b:b:b:b:b:b:b,c:c:c:c)}
%$$
%
%It is now obvious that, performing twice this operation,
%one gets back to the original partition $\lambda$, that is $c_n ^2 (\lambda )=\lambda$,
%as desired.\cvd
\cvd

As already mentioned, as a special case of the above proposition we have that the set of
all integer partitions is the set of fixed points of the map $c_1 ^2$ (where $c_1$ is the
conjugation map), because the condition of the proposition becomes empty in this case.
Another consequence is recorded in the following corollary,
which shows the role of Schr\"oder partitions in this context.

\begin{cor} The set of Schr\"oder partitions is the set of fixed points of the map
$c_2 ^2$.
\end{cor}

\emph{Proof.}\quad Just observe that, setting $n=2$ in the previous proposition,
the requirement in order to have $c_2 ^2 (\lambda )=\lambda$ is that
there is at most one part of $\lambda$ between two consecutive even numbers,
which means precisely that odd parts have to be simple.
\cvd

The set $\mathbf{Sch}$ of all Schr\"oder shapes can be naturally endowed with a poset structure,
by declaring $\lambda \leq \mu$ whenever the set of cells of the shape $\lambda$
is a subset of the set of cells of the shape $\mu$, provided that we draw the two shapes
in such a way that their top left cells coincide. This is equivalently (and perhaps more formally)
expressed in terms of Schr\"oder partitions: if $\lambda =(\lambda_1 ,\ldots ,\lambda_h )$ and
$\mu =(\mu_1 ,\ldots ,\mu_k )$, then $\lambda \leq \mu$ when $h\leq k$ and, for all
$i\leq h$, $\lambda_i \leq \mu_i$. Therefore the poset \Sch of Schr\"oder shapes is actually
a subposet of the Young lattice. However, it seems not at all a trivial one; notice, in particular,
that an interval of the Young lattice whose endpoints are Schr\"oder partitions
does not contain only Schr\"oder partitions (apart from very simple cases). In general,
it appears to be very hard (if not impossible) to infer nontrivial properties
of the Schr\"oder poset from properties of the Young lattice. The rest of this section is devoted to
developing some elements of the theory of the Schr\"oder poset along the lines suggested by
the classical theory of its more noble relative, the Young lattice.

\bigskip

One of the most fundamental properties the Schr\"oder poset shares with the Young lattice
is the fact that it is a distributive lattice. We will obtain this result as a consequence
of a more general one, which is of independent interest and can be seen as a slight generalization
of Lemma 2.1 in \cite{A}\footnote{We also notice that the Schr\"oder poset is a partition ideal of order 1,
in the terminology introduced in \cite{A}.}.

\begin{teor}
Given a function $f:\mathbf{N}\rightarrow \mathbf{N}\cup \{ \infty \}$,
denote with $\mathbf{P}_f$ the set of integer partitions in which part $i$ appears at most $f(i)$ times.
Then $\mathbf{P}_f$ is a distributive sublattice of the Young lattice (with partwise join and meet).
\end{teor}

\emph{Proof.}\quad Since every sublattice of a distributive lattice is distributive,
it will be enough to show that $\mathbf{P}_f$ is a sublattice of the Young lattice.

Given two partitions $\lambda, \mu \in \mathbf{P}_f$, their join in the Young lattice is
the partition $\lambda \vee \mu$ whose $i$-th part is the maximum between
$\lambda_i$ and $\mu_i$, for all $i$. We will now show that $\lambda \vee \mu$ is in $\mathbf{P}_f$.

By contradiction, suppose that part $i$ appears more than $f(i)$ ($\neq \infty$) times in $\lambda \vee \mu$,
and denote with $(\lambda \vee \mu )_{t}=\cdots =(\lambda \vee \mu)_{t+j}=i$ all such parts in $\lambda \vee \mu$.
Since $\lambda, \mu \in \mathbf{P}_f$, it cannot happen that $(\lambda \vee \mu)_{t+s}=\lambda_{t+s}$,
for all $s=0,\ldots j$, and the same holds with $\lambda$ replaced by $\mu$.
In other words, there exist two indices $k,r$, with $0<k\leq r<j$, such that (without loss of generality)
$\lambda_t =\lambda_{t+1}=\cdots =\lambda_{t+r}=i>\lambda_{t+r+1}$ and $\mu_{t+k-1}>\mu_{t+k} =\mu_{t+k+1}=\cdots =\mu_{t+j}=i$.
But this would imply, in particular, that $i=\lambda_{t+k-1}\geq \mu_{t+k-1}>\mu_{t+k}=i$, which is plainly impossible.
We can thus conclude that $\lambda \vee \mu$ is in $\mathbf{P}_f$.

Using a completely similar argument one can also show that
the meet of two partitions belonging to $\mathbf{P}_f$ in the Young lattice is again a partition of $\mathbf{P}_f$,
thus completing the proof.\cvd

\begin{cor}
The Schr\"oder poset \Sch is a distributive lattice.
\end{cor}

\emph{Proof.}\quad Just apply the previous theorem with $f$ defined by setting
$f(n)=\infty$ when $n$ is even and $f(n)=1$ when $n$ is odd.\cvd

\bigskip

The Young lattice is the prototypical example of a differential poset. Following \cite{St},
an \emph{r-differential poset} (for some positive integer $r$), is a locally finite,
ranked poset $\mathcal{P}$ having a minimum and such that:
\begin{itemize}
  \item for any two distinct elements $x,y$ of $\mathcal{P}$, if there are exactly $k$ elements
  covered by both $x$ and $y$, then there are exactly $k$ elements which cover both $x$ and $y$;
  \item if $x$ covers exactly $k$ elements, then $x$ is covered by exactly $k+r$ elements.
\end{itemize}

The Young lattice is a 1-differential poset. More specifically, it is the unique
1-differential distributive lattice. Thus, it is clear that \Sch is not a 1-differential poset,
since we have proved right now that it is a distributive lattice (and it is clearly not isomorphic
to the Young lattice). However, it belongs to
a wider class of posets which we believe to be an interesting generalization of differential posets.

\bigskip

Let $\varphi$ be a map sending a positive integer $k$ to an interval $\varphi (k)$ of positive integers.
We say that a poset $\mathcal{P}$ is a \emph{$\varphi$-differential poset} when
it is an infinite, locally finite, ranked poset with a minimum such that:
\begin{enumerate}
  \item for any two distinct elements $x,y$ of $\mathcal{P}$, if there are exactly $k$ elements
  covered by both $x$ and $y$, then there are exactly $k$ elements which cover both $x$ and $y$;
  \item if $x$ covers exactly $k$ elements, then $x$ is covered by $l$ elements,
  for some $l\in \varphi(k)$.
\end{enumerate}

When there exists a positive integer $r$ such that $\varphi (k)=\{
k+r\}$, for all $k$, a $\varphi$-differential poset is just an
$r$-differential poset.

\bigskip

The next proposition shows that \Sch is indeed a $\varphi$-differential distributive lattice,
for a suitable $\varphi$.

\begin{prop}\label{differential} Let $\lambda$ be a Schr\"oder partition covering $k$ Schr\"oder partitions in \Sch.
Then $\lambda$ is covered by $l$ Schr\"oder partitions, with $\lceil \frac{k+1}{2}\rceil \leq l\leq 2k$.
\end{prop}

\emph{Proof}.\quad Given $\lambda$ in \Sch, we denote with
$\uparrow \! \! \lambda$ the number of elements of \Sch covering
$\lambda$ and with $\downarrow \! \! \lambda$ the number of
elements of \Sch which are covered by $\lambda$. From the
hypothesis we have that $\downarrow \! \! \lambda =k$.

In the rest of the proof we slightly modify our notation for
partitions. Namely, we will add to each partition a smallest part
equal to 0. So, for instance, we will write $\lambda =(\lambda_1
,\lambda_2 ,\ldots ,\lambda_r ,\lambda_{r+1})$, with
$\lambda_{r+1}=0$. A part $\lambda_i$ of $\lambda$ will be called
\emph{up-free} whenever either
\begin{itemize}
\item it is odd, or
\item it is even and
$\lambda_{i-1}\neq \lambda_i ,\lambda_i+1$.
\end{itemize}
Similarly, it will be called \emph{down-free} whenever either
\begin{itemize}
\item it is odd, or
\item it is even and $\lambda_{i+1}\neq \lambda_i ,\lambda_i -1$.
\end{itemize}
In particular, each odd part of $\lambda$ is both up-free and down-free
(for this reason, we will sometimes refer to odd parts as \emph{trivial}
up-free (or down-free) parts). Observe that, concerning the special
(even) part $\lambda_{r+1}=0$, it is never down-free by
convention, whereas it is assumed to be up-free when $\lambda_r
\neq 1$.

In order to determine $\uparrow \! \! \lambda$, we observe that we
have to find those parts of $\lambda$ to which we can add 1
without losing the property of being a Schr\"oder partition. These
are precisely all up-free parts. Similarly, $\downarrow \! \!
\lambda$ is given by the number of down-free parts.

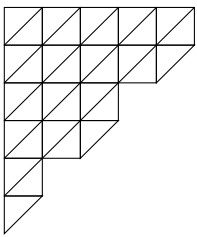
\begin{figure}
\begin{center}
\begin{tikzpicture}[scale=0.5]
%\shape{4,3,3,1}
\shape{5,4,3,2,1}

\draw (4,-1) -- (5,-1) -- (4,-2) -- cycle;

\draw (2,-3) -- (3,-3) -- (2,-4) -- cycle;

\draw (0,-5) -- (1,-5) -- (0,-6) -- cycle;
\end{tikzpicture}
\end{center}
\caption{A Schr\"oder shape which maximizes $\uparrow \! \!
\lambda$.\label{maximize}}
\end{figure}

If we want to maximize $\uparrow \! \! \lambda$, it is then clear
that we have to choose a Schr\"oder partition $\lambda$ having
many nontrivial up-free parts and few nontrivial down-free parts
(odd parts are irrelevant). Observe moreover that we can restrict
ourselves to the case of $\lambda$ having all distinct parts,
since several repeated (even) parts is equivalent to having
only one part of the same cardinality (all parts except for the
top one cannot be modified). Concerning the
greatest part of $\lambda$, $\lambda_1$, we notice that it has to
be even, otherwise it would be down-free. Moreover, in order to
have few partitions immediately below $\lambda$, we should try to
make $\lambda_1 =2n$ not down-free. To do this, just choose
$\lambda_2 =\lambda_1 -1=2n-1$ (recall that we are assuming
$\lambda$ to have all distinct parts). Observe that, in this way,
$\lambda_2$ is odd (and so trivially down-free), however it is not
difficult to realize that any other choice of $\lambda_2$ would
have produced a down-free part (an even part strictly larger than
another even part is certainly down-free). Now, concerning
$\lambda_3$, we wish it to be up-free but not down-free. The first
condition is fulfilled if and only if $\lambda_3 \neq \lambda_2
-1=2n-2$; for the second condition, we must choose $\lambda_3$
even and such that $\lambda_4 =\lambda_3 -1$. Without loss of
generality, we can set $\lambda_3 =2n-4$, so that $\lambda_4
=2n-5$. We can now argue in a completely analogous way for all the
remaining parts of $\lambda$, until we have $\downarrow \! \!
\lambda =k$ (notice that $n$ has to be large enough to reach this
goal). In the end, we obtain that a partition $\lambda$ which
maximizes $\uparrow \! \! \lambda$ has odd-indexed parts
$\lambda_{2i+1}=2n-4i$ and even-indexed parts
$\lambda_{2i}=2n-4i-1$ (see Figure \ref{maximize} for an example).

%For instance, one such partition is $(14,13,10,9,6,5)$, which covers 3
%Schr\"oder partitions (namely $(14,12,10,9,6,5)$,
%$(14,13,10,8,6,5)$ and $(14,13,10,9,6,4)$ and is covered by 7
%Schr\"oder paritions (namely $(15,13,10,9,6,5)$,
%$(14,14,10,9,6,5)$, $(14,13,11,9,6,5)$, $(14,13,10,8,6,5)$,
%$(14,13,10,9,7,5)$, $(14,13,10,9,6,6)$ and $(14,13,10,9,6,5,1)$).
A direct computation then shows that, in the best possible cases
(which occur when the smallest part of $\lambda$ is $1$ or $2$),
we get $\uparrow \! \! \lambda =2k$, as desired.

A similar approach allows also to determine a lower bound for
$\uparrow \! \! \lambda$. The only difference with the previous
arguments is that now we would like to have a partition $\lambda$
having many down-free parts and few up-free parts. It turns out
that the role of odd and even parts are somehow swapped in the
above arguments. Specifically, it can be shown that the largest
part $\lambda_1$ of $\lambda$ has to be odd, and that $\lambda_2
=\lambda_1 -1$. At the end, we obtain a partition having
odd-indexed parts $\lambda_{2i+1}=2n-4i-1$ and even-indexed parts
$\lambda_{2i}=2n-4i-2$. Similarly as before, a direct computation
shows that, when the smallest part of $\lambda$ is $1$ or $2$, we
get the desired lower bound. The task of providing all the details
is then left to the reader.\cvd

\section{An RSK-like algorithm for Schr\"oder tableaux}\label{algorithm}

From the algorithmic point of view, the main application of Young tableaux is in the context of the RSK algorithm.
This algorithm, named after Robinson, Schensted and Knuth, takes as input
a word (on the alphabet of positive integers) of length $n$
and produces in output two semistandard Young tableaux with $n$ cells having the same shape.
For what concerns us, we will deal with a special case of the RSK algorithm,
often referred to as \emph{Robinson-Schensted correspondence} (briefly, RS correspondence), in which the input is
a permutation of length $n$ and the output is given by a pair of standard Young tableaux.
A brief description of such an algorithm is given below (Algorithm \ref{RS},
where $\pi =\pi_1 \pi_2 \cdots \pi_n$ is a generic permutation of length $n$).

\begin{algorithm}[!h]
%\KwIn{sane people and photos} \KwOut{bebo addicts}
$P:={\small \young(\pione )}$\; $Q:={\small \young(1)}$\;
\For{\textnormal{$k$ from 2 to $n$}}
    {
    $\alpha :=\pi_k$\;
    \For{$i\geq 1$}
        {
        \eIf
            {
            \textnormal{$\alpha$ is bigger than all elements in the $i$-th row of $P$}
            }
            {
            append a cell with $\pi_k$ inside at the end of the $i$-th row of $P$\;
            append the cell ${\tiny \young(k)}$ at the end of the $i$-th row of $Q$\;
            break\;
            }
            {
            write $\alpha$\ in the cell of the $i$-th row containing the smallest element $\beta$ bigger than $\alpha$\;
            $\alpha :=\beta$\;
            }
        }
    }
%$S:=\emptyset$\; $Q:=E$\; \While{$Q\neq \emptyset$}{find $m\in Q$
%having maximum weight\; $Q:=Q\setminus \{ m\}$\; \If{$S\cup \{ m\}
%\in \mathcal{F}$}{$S:=S\cup \{ m\}$\;}} \Return{$S$}\;
\caption{RS($\pi$)}\label{RS}
\end{algorithm}

The RSK algorithm is extensively described in the literature. For
instance, the interested reader can find a modern and elegant
presentation of it in \cite{Be}. Among other things, one of the
most beautiful properties of the RS correspondence is that it
establishes a bijection between permutations of length $n$ and
pairs of standard Young tableaux with $n$ cells having the same
shape. This fact bears important enumerative consequences,
as well as strictly algebraic ones. For a given permutation $\pi$, the tableaux of the pair$(P,Q)$
returned by the RS algorithm are usually referred to as the
\emph{insertion tableau} (the tableau $P$) and the \emph{recording
tableau} (the tableau $Q$). As a consequence, we have the
following nice result, which can again be found in \cite{Be}.

\begin{teor} Denote with $f^\lambda$ the number of standard Young tableaux of shape $\lambda$. Then we have:
$$
n!=\sum_{\lambda \vdash n}(f^\lambda)^2 .
$$
\end{teor}

A \emph{standard Schr\"oder tableau} (from now on, simply \emph{Schr\"oder tableau})
with $n$ cells is a Schr\"oder shape whose cells are filled in
with positive integers from $1$ to $n$ in such a way that entries
in each row and each column are (strictly) increasing.

We propose here a natural analog of the RS algorithm for Schr\"oder tableaux. The main difference (which is due
to the specific underlying shape of a Schr\"oder tableaux) lies in the fact that there are two distinct ways of
managing the insertion of a new element in the tableau, depending on whether the cell it should be inserted in
is an upper triangle or a lower triangle. As a consequence, our algorithm does not establish a bijection
between permutations and pairs of Schr\"oder tabealux; nevertheless, due to the strict analogy
with the RS correspondence, we believe that it is very likely to have interesting combinatorial properties.
A description of our algorithm is given below (Algorithm \ref{Sch}, where $\pi$ is as in Algorithm \ref{RS}).

\begin{algorithm}[!h]
%\KwIn{sane people and photos} \KwOut{bebo addicts}
$P:=$ the 1-cell Schr\"oder tableau with $\pi_1$ written in the cell\;
$Q:=$ the 1-cell Schr\"oder tableau with $1$ written in the cell\;
\For{\textnormal{$k$ from 2 to $n$}}
    {
    $\alpha :=\pi_k$\;
    \For{$i\geq 1$}
        {
        \eIf
            {
            \textnormal{$\alpha$ is bigger than all elements in the $i$-th row of $P$}
            }
            {
            append a cell (either an upper or a lower triangle) with $\pi_k$ inside at the end of the $i$-th row of $P$\;
            append a cell (either an upper or a lower triangle) with $k$ inside at the end of the $i$-th row of $Q$\;
            break\;
            }
            {let $A$ be the cell of the $i$-th row containing the smallest element bigger than $\alpha$\;
            \eIf
                {
                \textnormal{$A$ is an upper triangle}
                }
                {
                if $A$ has a twinned lower triangle immediately below it, set $\beta :=$ content of such lower triangle\;
                move the content of $A$ to the lower triangle immediately below $A$, possibly by creating such a triangle if it does not exist\;
                write $\alpha$ in $A$\;
                if $\beta$ exists, set $\alpha :=\beta$, else break\;
                }
                {
                $\beta :=$ content of $A$\;
                write $\alpha$ in $A$\;
                $\alpha :=\beta$\;
                }

            }
        }
    }
%$S:=\emptyset$\; $Q:=E$\; \While{$Q\neq \emptyset$}{find $m\in Q$
%having maximum weight\; $Q:=Q\setminus \{ m\}$\; \If{$S\cup \{ m\}
%\in \mathcal{F}$}{$S:=S\cup \{ m\}$\;}} \Return{$S$}\;
\caption{Sch($\pi$)}\label{Sch}
\end{algorithm}

\bigskip

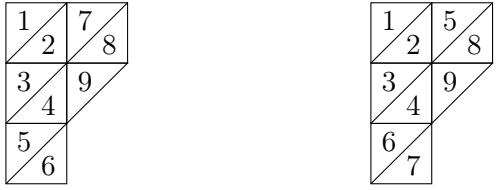
\begin{figure}
\begin{center}
\begin{tikzpicture}[scale=0.8]
\shape{2,1,1}

\draw (1,-1) -- (2,-1) -- (1,-2) -- cycle;

\node[] (1) at (0.3,-0.3) {1};
\node[] (2) at (0.7,-0.7) {2};
\node[] (7) at (1.3,-0.3) {7};
\node[] (8) at (1.7,-0.7) {8};
\node[] (3) at (0.3,-1.3) {3};
\node[] (4) at (0.7,-1.7) {4};
\node[] (9) at (1.3,-1.3) {9};
\node[] (5) at (0.3,-2.3) {5};
\node[] (6) at (0.7,-2.7) {6};

\secondshape{2,1,1}

\draw (7,-1) -- (8,-1) -- (7,-2) -- cycle;

\node[] (1) at (6.3,-0.3) {1};
\node[] (2) at (6.7,-0.7) {2};
\node[] (5) at (7.3,-0.3) {5};
\node[] (8) at (7.7,-0.7) {8};
\node[] (3) at (6.3,-1.3) {3};
\node[] (4) at (6.7,-1.7) {4};
\node[] (9) at (7.3,-1.3) {9};
\node[] (6) at (6.3,-2.3) {6};
\node[] (7) at (6.7,-2.7) {7};

\end{tikzpicture}
\end{center}
\caption{How our RS-like algorithm works.}\label{RSKlike}
\end{figure}

\emph{Example.}\quad Consider the permutation $\pi =465193287$. The pair $(P,Q)$ of Schr\"oder tableaux produced
by applying the algorithm Sch to $\pi$ is illustrated in Figure \ref{RSKlike}.

\bigskip

In this section we aim at starting the investigation of the combinatorial properties of this
RS-analog. More specifically, we will address the following problems: given a Schr\"oder shape
$P$, can we characterize those permutations having a Schr\"oder tableau of shape $P$ as their insertion tableau?
How many of them are there? This problem seems to be quite difficult in its full generality; here we will deal
with very few simple cases, for which we can provide complete answers.

\subsection{Permutations with given Schr\"oder insertion shape: some cases}

In the present subsection we collect some starting results concerning permutations whose Schr\"oder insertion tableaux
have simple shapes. Some of the results will be stated without proof, since they have already been proved in \cite{Fe}.

The first case we investigate is that of a Schr\"oder shape consisting of a single row
(which can terminate either with an upper or a lower triangle). To state our result
we first need to recall a classical definition.

Given a permutation $\pi =\pi_1 \cdots \pi_n$, we say that $\pi_i$ is a
\emph{left-to-right maximum} (or, briefly, \emph{LR maximum}) whenever
$\pi_i =\max (\pi_1 ,\ldots ,\pi_i)$.

\begin{prop}
Let $\pi =\pi_1 \cdots \pi_n$ be a permutation of length $n$.
The Schr\"oder insertion tableau of $\pi$ has a single row if and only if,
for all $i\leq n$:
\begin{enumerate}
\item if $i$ is odd, then $\pi_i$ is a LR maximum of $\pi$;
\item if $i$ is even, then $\pi_i$ is a LR maximum of the permutation
obtained from $\pi$ by removing $\pi_{i-1}$
(and suitably renaming the remaining elements).
\end{enumerate}
\end{prop}

%\emph{Proof.}\quad Suppose we are inserting $\pi_i$ in the insertion tableau $P$,
%which is assumed to consist of a single row. If $i$ is odd, then the last cell of $P$
%is a lower triangle; in order not to create new rows, $\pi_i$ has necessarily to be a LR maximum.
%On the other hand, if $i$ is even, then the last cell of $P$ is an upper triangle;
%in this case, $\pi_i$ can be inserted in $P$ in two ways: either $\pi_i$ is a LR maximum,
%and so it is simply appended at the end of the unique row of $P$, or $\pi_i$ is greater than
%all previous elements of $\pi$ but $\pi_{i-1}$, hence $\pi_i$ is inserted in the cell
%containing $\pi_{i-1}$ (which is the last cell of the unique row of $P$, and so it is
%an upper triangle) and a new cell (a lower triangle) containing $\pi_{i-1}$ is added
%at the end of the unique row of $P$.
%
%Conversely, it is easy (and so left to the reader) to check that a permutation
%satisfying conditions 1 and 2 in the statement of the present proposition
%must have a Schr\"oder insertion tableau consisting of a single row.\cvd

The permutations $\pi$ of length $n$ whose Schr\"oder insertion
tableau have a single row can therefore be simply characterized as
follows: for all $i$, $\{ \pi_{2i+1},\pi_{2i+2}\} =\{
2i+1,2i+2\}$. As a consequence of this fact, a formula for the number of such permutations
follows immediately.

\begin{prop}
The set of permutations of length $n$ whose Schr\"oder insertion tableau consists of a single row
has cardinality $2^{\lfloor \frac{n}{2}\rfloor}$.
\end{prop}

The second case we consider is the natural counterpart of the previous one, that is
Schr\"oder shapes having a single column. Despite the similarities with the previous case,
it turns out that the set of permutations having Schr\"oder insertion tableau of this form
can be nicely described in terms of \emph{pattern avoidance}.

Given two permutations $\sigma$ and $\tau =\tau_1 \cdots \tau_n$
(of length $k$ and $n$ respectively, with $k\leq n$), we say that
there is an \emph{occurrence} of $\sigma$ in $\tau$ when there
exists indices $i_1 <i_2 <\cdots <i_k$ such that
$\tau_{i_1}\tau_{i_2}\cdots \tau_{i_k}$ is order isomorphic to
$\sigma$. When there is an occurrence of $\sigma$ in $\tau$, we
also say that $\tau$ \emph{contains the pattern} $\sigma$. When
$\tau$ does not contain $\sigma$, we say that $\tau$ \emph{avoids
the pattern} $\sigma$. The set of all permutations of length $n$
avoiding a given pattern $\sigma$ is denoted with $Av_n (\sigma
)$. Some useful references for the combinatorics of patterns in
permutations are \cite{Bo} and \cite{Ki}, whereas similar notions
of patterns in set partitions and in compositions and words are studied
in \cite{M} and \cite{HM}, respectively.

\begin{prop}\label{1column}
Let $\pi =\pi_1 \cdots \pi_n$ be a permutation of length $n$.
The Schr\"oder insertion tableau of $\pi$ has a single column if and only if
$\pi \in Av_n (123,213)$.
\end{prop}

%\emph{Proof.}\quad An argument similar to that of the preceding proposition shows that
%the Schr\"oder insertion tableau of $\pi$ has a single column if and only if, for all $i\leq n$,
%$\pi_i <\min (\{\pi_1 ,\ldots ,\pi_{i-1}\} \setminus \min \{ \pi_1 ,\ldots ,\pi_{i-1}\} )$
%(i.e., $\pi_i$ is smaller than the second minimum of set of all previous elements).
%Thus $\pi$ can be factored into subpermutations (made of consecutive elements of $\pi$),
%say $\pi =\tilde{\pi_1}\cdots \tilde{\pi_r}$, in such a way that each factor $\tilde{\pi_i}$
%is isomorphic to a permutation of the form $1t(t-1)\cdots 32$ (for some $t$) and each element
%of $\tilde{\pi_i}$ is greater than each element of $\tilde{\pi_{i+1}}$ (for all $i$).
%In the language of permutation patterns, this is usually expressed by saying that
%$\pi$ is a \emph{skew sum} of permutations of the form $1t(t-1)\cdots 32$.
%It is now a known fact (see, for instance, \cite{AA}) that such permutations are precisely those
%avoiding the two patterns $123$ and $213$.\cvd

Many classes of permutations avoiding a given set of patterns have been enumerated.
The above one is among them, see \cite{SiSc}.

\begin{prop}
The set of permutations of length $n$ whose Schr\"oder insertion tableau consists of a single column
has cardinality $2^{n-1}$.
\end{prop}

We close this section by illustrating one more case which is, in some sense, a generalization of both the cases described above.
Namely, we consider the case of what can be called \emph{Schr\"oder hooks},
that is Schr\"oder shapes having at most one row and one column with more than one upper triangular cell.
Since this case is considerably more difficult than the previous ones, we need some preparation and our results will be less elegant.
Nevertheless, the strategy employed reveals some features of our algorithms that are interesting in themselves.

Given a permutation $\pi =\pi_1 \pi_2 \cdots \pi_n$ (that will be fixed until the end of this section), $\pi_i$ is a \emph{quasi-left-to-right minimum} (briefly, a QLTR minimum) of $\pi$ when either $i=1,2$ or $\pi_i <\min (\{ \pi_1 ,\ldots ,\pi_{i-1}\} \setminus \{ \min (\pi_1 ,\ldots ,\pi_{i-1})\})$.
In other words, a QLTR minimum of $\pi$ is an element of $\pi$ which is smaller than the second smallest element of $\pi$ preceding it.
For instance, in the permutation
$\rho =\underline{4}\,\underline{1}\,7\,\underline{2}\,9\,3\,6\,5\,8$,
the QLTR minima are the underlined elements, i.e. $4,1,2$.

Given a natural number $i$, the \emph{$i$-th QLTR sequence of $\pi$} is recursively defined as follows:
\begin{itemize}
\item when $i=1$, it is the sequence of the QLTR minima of $\pi$;
\item when $i>1$, it is the sequence of the QLTR minima of the permutation obtained from $\pi$ by deleting the elements of its $j$-th QLTR sequence,
for all $j<i$.
\end{itemize}

The $i$-th QLTR sequence of a permutation can be interpreted as a permutation as well, by simply replacing its $k$-th smallest element with $k$.
When no confusion is likely to arise, we will call ``$i$-th QLTR sequence" both the sequence and the associated permutation.

Every permutation $\pi$ can be written as the shuffle of its $QLTR$ sequences.
Considering again the permutation $\rho$ above, such a decomposition is the following:
$\rho =\underline{4}\,\underline{1}\,\widetilde{7}\,\underline{2}\,\widetilde{9}\,\widetilde{3}\,\widetilde{6}\,\widetilde{5}\,\widehat{8}$.

\begin{lemma}
The $i$-th QLTR sequence of a permutation $\pi$ avoids 123 and 213, for all $i$.
\end{lemma}

\emph{Proof.}\quad Just observe that, for each $i$,
the $i$-th sequence of $\pi$ consists of thoes elements which enter the Schr\"oder insertion tableau of $\pi$ in column $i$.
So, using an argument completely analogous to that of Proposition \ref{1column} (see \cite{Fe}), we get the thesis.\cvd

We now label the elements of $\pi$ by recording the column in which they enter the Schr\"oder insertion tableau.
We say that the \emph{c-label} of $\pi_i$ is $c_j$ when it enters the Schr\"oder insertion tableau of $\pi$ in column $j$.
The \emph{c-word} of $\pi$ is then the word obtained from $\pi$ by replacing each $\pi_i$ with its c-label.
Moreover, the \emph{bumping word} of $\pi$ is obtained from by c-word by deleting the two occurrences of $c_j$
corresponding to the two smallest elements whose c-label is $c_j$, for every $j$,
and the \emph{bumping sequence} of $\pi$ is the sequence of elements of $\pi$ corresponding to its bumping word.
Therefore the c-word of our running example $\rho$ is $c_1 ^2 c_2 c_1 c_2 ^4 c_3$, its bumping word is $c_1 c_2 ^3$
and its bumping sequence is $4\,7\,9\,6$.

We come finally to our last definition. The \emph{ordered bumping sequence} is obtained from the bumping sequence of $\pi$
by rearranging in decreasing order the element having the same c-label (and keeping their relative positions).
So the ordered bumping sequence of $\rho$ is $4\,9\,7\,6$.

We are now ready to state our main result on Schr\"oder hooks.

\begin{prop}
The Schr\"oder insertion tableau of $\pi$ is a Schr\"oder hook if and only if the ordered bumping sequence of $\pi$ avoids 123 and 213.
\end{prop}

\emph{Proof.}\quad Saying that the Schr\"oder insertion tableau of $\pi$ is a Schr\"oder hook is equivalent to saying that
every element of $\pi$ which is bumped down from the first to the second row always goes to the first column.
Now observe that, for every $j$, the first two elements of $\pi$ having c-label $j$ are inserted into the $j$-th cell of the first row
without causing any element to be bumped down.
On the other hand, all the successive elements that are inserted into the same cell bump down an element having the same c-label.
This means that the bumping word of $\pi$ records the c-labels of the elements that are successively bumped down into the second row.
However, the bumping sequence \emph{is not} the sequence of the bumped elements.
Instead, each element $\pi_i$ of the bumping sequence having c-label $j$ bumps down
the second smallest element among those preceding $\pi_i$ and having c-label $j$.
In fact, when $\pi_i$ is inserted, all the previous elements having c-label $j$ have been bumped down, except for the two smallest ones,
and of course $\pi_i$ bumps down the largest of the two.
As a consequence, we observe that the elements of the $j$-th QLTR sequence of $\pi$ are bumped down in decreasing order,
with the two smallest ones which remain in the first row.
Therefore, summing up the above considerations, we have that the $k$-th element of $\pi$ which is bumped down from the first row
has c-label equal to the $k$-th letter of the bumping word;
moreover, since the set of elements having the same c-label are bumped down in decreasing order,
the ordered bumping sequence is precisely the sequence of the elements that are bumped down in the correct order.
Now, it is clear that the Schr\"oder insertion tableau of $\pi$ is a Schr\"oder hook if and only if
all the elements of $\pi$ that are bumped down from the first row are placed into the first column,
and this happens if and only if the ordered bumping sequence represents a permutation whose insertion Schr\"oder tableau has a single column.
From Proposition \ref{1column}, we know that this happens precisely when the ordered bumping sequence avoids both 123 and 213, as desired.\cvd

\section{An alternative view of Schr\"oder tableaux}\label{alternative}

Following our treatment, Schr\"oder tableaux can be interpreted as
upper saturated chains in Schr\"oder lattices (where upper means
that the maximum of the chain is the maximum of the lattice). Now we propose a
different description of Schr\"oder tableaux, relying on at least
two main ingredients: interval orders (which are a well known
class of posets) and a notion of \emph{weak pattern} for posets,
which is not entirely new in its own right, but which appears to
have never been considered from a strictly combinatorial point of view.

\subsection{Interval orders}\label{interval}

A poset $\mathcal{P}$ is called an \emph{interval order} when it
is isomorphic to a collection of intervals of the real line, with
partial order relation given as follows: for any two intervals
$I,J$, it is declared that $I<J$ whenever all elements of $I$ are
less than all elements of $J$. In other words, the interval $I$
lies completely on the left of $J$. For the purposes of the
present article, all intervals will be closed, and the minimum and
the maximum will be natural numbers. Notice that, under these hypotheses,
the set of all maxima and minima of the intervals of a given interval order
can be chosen to be an initial segment of the natural numbers.

The notion of interval order is now very classical, and was
introduced by Fishburn \cite{Fi}. Though the main motivation for the
introduction of such a concept came from social choice theory,
it soon revealed its intrinsic interest, especially from a
combinatorial point of view. To support this statement (and without
giving any detail), we only recall here the characterization of
interval orders as partially ordered sets avoiding the (induced) subposet
$\mathbf{2}+\mathbf{2}$, and the more recent enumeration of
finite interval orders \cite{BCDK}.

An immediate link between interval orders and Schr\"oder tableaux
is given by the fact that every Schr\"oder tableau can be associated
with a set of intervals. Given a Schr\"oder shape $\lambda$,
two cells $A$ and $B$ of $\lambda$ are called \emph{twin} when
they are adjacent and their union is a square. Equivalently,
two adjacent cells $A$ and $B$ are twin cells when their common
edge is a diagonal edge. Moreover, a (necessarily upper triangular)
cell of $\lambda$ is called \emph{lonely} when it is the last cell of an odd row.
Notice that the set of cells of $\lambda$ can be partitioned into twin cells
and lonely cells. Now, given a Schr\"oder tableau $S$ having $n$ cells,
consider the set of intervals $\mathcal{I}_S$ defined as follows:
$I=[a,b]\in \mathcal{I}_S$ when both $a$ and $b$ are fillings of
a pair of twin cells of $S$ or $a$ is the filling of a lonely cell and $b=n+1$.
For instance, for the Schr\"oder tableau $S$ on the right in Figure \ref{RSKlike},
which has 9 cells, we have $\mathcal{I}_S =\{ [1,2],[5,8],[3,4],[9,10],[6,7]\}$.
The benefit of endowing the set of intervals associated with a Schr\"oder tableau
with its interval order will be discussed in the next subsections.

\subsection{Weak patterns in posets and strong pattern avoidance}

The study of classes of posets which contain or avoid certain subposets is a major trend
in order theory. Classically, a poset $\mathscr{Q}$ \emph{contains} another poset $\mathscr{P}$
whenever $\mathscr{Q}$ has a subposet isomorphic to $\mathscr{P}$. Borrowing the terminology from permutations,
we could also say that $\mathscr{Q}$ contains the \emph{pattern} $\mathscr{P}$. On the other hand,
we say that $\mathscr{Q}$ \emph{avoids} $\mathscr{P}$ whenever $\mathscr{Q}$ does not contain $\mathscr{P}$.
The notion of pattern containment defines a partial order on the set $\mathfrak{X}$ of all (finite) posets,
and we will write $\mathscr{P} \sqsubseteq \mathscr{Q}$ to mean that $\mathscr{P}$ is contained in $\mathscr{Q}$.
Instead, the class of all finite posets avoiding a given poset $\mathscr{P}$ will be denoted
$Av(\mathscr{P})$.

Here the use of the word ``subposet'' might be controversial. Technically speaking,
what we have called ``subposet'' is sometimes called ``induced subposet''. Formally, we say that
$\mathscr{P}$ is an \emph{induced subposet} of $\mathscr{Q}$ when there is an injective function
$f:\mathscr{P} \rightarrow \mathscr{Q}$ which is both order-preserving and order-reflecting:
for all $x,y$, $x\leq y$ in $\mathscr{P}$ if and only if $f(x)\leq f(y)$ in $\mathscr{Q}$.
Loosely speaking, this means that $\mathscr{Q}$ contains an isomorphic copy of $\mathscr{P}$.
In what follows we will fully adhere to such terminology: $\mathscr{P} \sqsubseteq \mathscr{Q}$
thus means that $\mathscr{Q}$ has an induced subposet isomorphic to $\mathscr{P}$.

What is useful for us is however a weaker version of the above notion of pattern.
We say that $\mathscr{P}$ is \emph{weakly contained} in $\mathscr{Q}$ (or that
$\mathscr{P}$ is a \emph{weak pattern} of $\mathscr{Q}$) when there exists
an injective order-preserving function $f:\mathscr{P}\rightarrow \mathscr{Q}$.
This can be also expressed by saying that $\mathscr{P}$ is a (not necessarily induced) subposet
of $\mathscr{Q}$. It is clear that a pattern is also a weak pattern.
On the other hand, we say that $\mathscr{Q}$ \emph{strongly avoids} $\mathscr{P}$ whenever
$\mathscr{Q}$ does not weakly contain $\mathscr{P}$. This is also expressed by writing
$\mathscr{Q}\in SAv(\mathscr{P})$.

\begin{figure}
\begin{center}
\begin{tikzpicture}[scale=0.4]

%\draw (1,-1) -- (2,-1) -- (1,-2) -- cycle;

\draw (5,1) [fill] circle (.1);
\draw (6,1) [fill] circle (.1);
\draw (7,1) [fill] circle (.1);
\draw (6,1) circle (1.5);

\draw (5.5,6) [fill] circle (.1);
\draw (6.5,6) [fill] circle (.1);
\draw (5.5,7) [fill] circle (.1);
\draw[very thin] (5.5,6) -- (5.5,7);
\draw (6,6.3) circle (1.5);

\draw (1,10) [fill] circle (.1);
\draw (2,11) [fill] circle (.1);
\draw (3,10) [fill] circle (.1);
\draw[very thin] (1,10) -- (2,11);
\draw[very thin] (2,11) -- (3,10);
\draw (2,10.2) circle (1.5);

\draw (9,11) [fill] circle (.1);
\draw (10,10) [fill] circle (.1);
\draw (11,11) [fill] circle (.1);
\draw[very thin] (9,11) -- (10,10);
\draw[very thin] (10,10) -- (11,11);
\draw (10,10.2) circle (1.5);

\draw (6,13.2) [fill] circle (.1);
\draw (6,14.2) [fill] circle (.1);
\draw (6,15.2) [fill] circle (.1);
\draw[very thin] (6,13.2) -- (6,14.2);
\draw[very thin] (6,14.2) -- (6,15.2);
\draw (6,14.2) circle (1.5);

\draw[thick] (6,2.5) -- (6,4.8);
\draw[thick] (5.07,7.30) -- (3.31,9.38);
\draw[thick] (6.93,7.30) -- (8.69,9.38);
\draw[thick] (3.31,10.80) -- (5.07,12.88);
\draw[thick] (8.69,10.80) -- (6.93,12.88);

\draw (19,1) [fill] circle (.1);
\draw (19,3) [fill] circle (.1);
\draw (17,5) [fill] circle (.1);
\draw (19,5) [fill] circle (.1);
\draw (21,5) [fill] circle (.1);
\draw (15,7) [fill] circle (.1);
\draw (18,7) [fill] circle (.1);
\draw (20,7) [fill] circle (.1);
\draw (23,7) [fill] circle (.1);
\draw (17,9) [fill] circle (.1);
\draw (19,9) [fill] circle (.1);
\draw (21,9) [fill] circle (.1);
\draw (17,11) [fill] circle (.1);
\draw (19,11) [fill] circle (.1);
\draw (21,11) [fill] circle (.1);
\draw (19,13) [fill] circle (.1);

\draw (19,1) -- (19,3);
\draw (19,3) -- (17,5);
\draw (19,3) -- (19,5);
\draw (19,3) -- (21,5);
\draw (17,5) -- (15,7);
\draw (17,5) -- (18,7);
\draw (17,5) -- (20,7);
\draw (19,5) -- (20,7);
\draw (21,5) -- (18,7);
\draw (21,5) -- (20,7);
\draw (21,5) -- (23,7);
\draw (15,7) -- (17,9);
\draw (18,7) -- (17,9);
\draw (18,7) -- (21,9);
\draw (20,7) -- (17,9);
\draw (20,7) -- (19,9);
\draw (20,7) -- (21,9);
\draw (23,7) -- (21,9);
\draw (17,9) -- (17,11);
\draw (17,9) -- (19,11);
\draw (19,9) -- (17,11);
\draw (19,9) -- (21,11);
\draw (21,9) -- (19,11);
\draw (21,9) -- (21,11);
\draw (17,11) -- (19,13);
\draw (19,11) -- (19,13);
\draw (21,11) -- (19,13);

\end{tikzpicture}
\end{center}
\caption{Hasse diagrams of $\mathfrak{X}_3$ (with explicit representation of each element)
and $\mathfrak{X}_4$.}\label{weak_inclusion}
\end{figure}
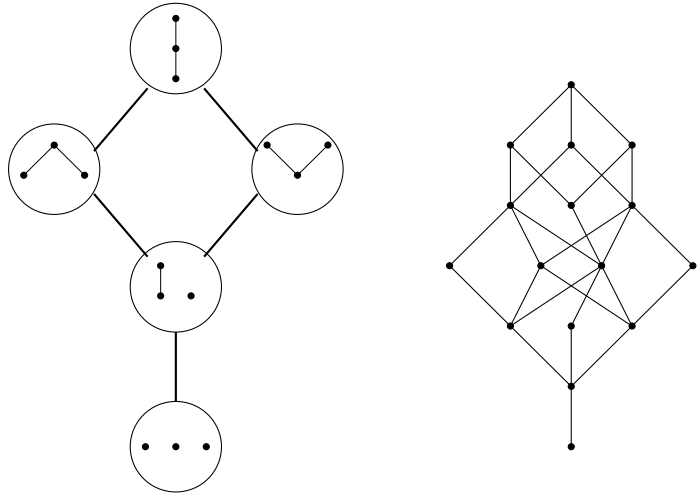

The partial order relation (on the set $\mathfrak{X}$ of all finite posets) defined by weak containment
will be denoted $\leq$. This notion of poset containment is not entirely new. Typically,
it has been considered in the context of families of sets, rather than generic posets,
and many investigations in this field concern the study of finite families of sets which strongly avoid
one or more finite posets, often with special focus on extremal properties (see for instance
\cite{GL,KT}). Here, however, we consider this order relation from a
purely order-theoretic point of view, with the aim of initiating the investigation of the poset
$(\mathfrak{X},\leq )$. Some rather easy facts are the following:

\begin{itemize}

\item $(\mathfrak{X},\leq )$ has minimum, which is the empty poset, and does not have maximum.

\item Given $\mathscr{P},\mathscr{Q}\in \mathfrak{X}$, if $\mathscr{P}\leq \mathscr{Q}$ then
the ground set of $\mathscr{P}$ has at most as many elements as $\mathscr{Q}$.

\item $(\mathfrak{X},\leq )$ is a ranked poset, and the rank function is the sum of the
number of elements of the ground set and the number of order relations between them.

\item Denoting $\mathfrak{X}_n$ the set of all posets of
size\footnote{The size is the number of elements of the ground set.} $n$, the restriction of
$\leq$ to $\mathfrak{X}_n$ gives a poset with minimum (the discrete poset on $n$ elements) and
maximum (the chain having $n$ elements). Also, $\mathfrak{X}_n$ has exactly one atom,
which is the poset of size $n$ having a single covering relation. Notice that,
if we replace $\leq$ with $\sqsubseteq$,
the resulting poset structure on $\mathfrak{X}_n$ would be trivial (more precisely, discrete).
We claim that the study of the posets $(\mathfrak{X}_n ,\leq )$
might be a potentially very interesting field of research.
In Figure \ref{weak_inclusion} we illustrate the posets $(\mathfrak{X}_n ,\leq )$
for a couple of small values of $n$.

\end{itemize}

Main aim of the present section is to initiate the study of the notion of strong pattern avoidance
introduced above. Recall that $SAv(\mathscr{P})$ denotes the class of all posets strongly avoiding
$\mathscr{P}$. Moreover, $SAv_n (\mathscr{P})$ is the subset of $SAv(\mathscr{P})$ consisting of
the posets of size $n$. Finally, the above expression can be easily adapted to the case of strong avoidance
of several posets (just by listing all the posets which are required to be avoided).
Observe that, for given posets $\mathscr{P}_1 ,\mathscr{P}_2 ,\ldots ,\mathscr{P}_s$,
$SAv(\mathscr{P}_1 ,\mathscr{P}_2 ,\ldots ,\mathscr{P}_s )$ is a \emph{down-set} of $(\mathfrak{X},\leq )$,
which means that, if $\mathscr{Q}\in SAv(\mathscr{P}_1 ,\mathscr{P}_2 ,\ldots ,\mathscr{P}_s )$ and
$\mathscr{P}\leq \mathscr{Q}$, then $\mathscr{P}\in SAv(\mathscr{P}_1 ,\mathscr{P}_2 ,\ldots ,\mathscr{P}_s )$.
In what follows we will always deal with the strong avoidance of a single poset.
The generalization of a given statement to several posets is easy (when meaningful) and so it is left to the reader.

\begin{prop}\label{upset}
For every poset $\mathscr{P}\in \mathfrak{X}_n$, $SAv(\mathscr{P})=Av(\langle \mathscr{P}\rangle_n )$,
where $\langle \mathscr{P}\rangle_n$ is the up-set\footnote{In a poset $\mathscr{P}$, an up-set $F$
is a subset of $\mathscr{P}$ such that, if $x\in F$ and $y\geq x$, then $y\in F$.
For a given subset $A$ of $\mathscr{P}$, the up-set generated by $A$ is the smallest up-set of $\mathscr{P}$
containing $A$.} generated by $\mathscr{P}$ in $\mathfrak{X}_n$.
\end{prop}

\emph{Proof.}\quad The fact that $\mathscr{Q} \in SAv(\mathscr{P})$ is equivalent to the following:
there cannot be an induced subposet $\mathscr{R}$ of size $n$ of $\mathscr{Q}$ which contains $\mathscr{P}$ as a
subposet. In other words, this means that $\mathscr{Q}$ cannot contain any induced subposet $\mathscr{R}$ of size $n$
which is a refinement of $\mathscr{P}$, that is $\mathscr{Q}\in Av(\langle \mathscr{P}\rangle_n )$, as desired.\cvd

The above proposition, rather than telling that strong avoidance can just be expressed in terms of classical avoidance,
suggests us that the formalism of strong poset avoidance allows to express certain problems concerning classical avoidance
in a much simplified way: avoiding several posets is sometimes equivalent to strongly avoiding just one of them.

\bigskip

We next characterize some classes of posets strongly avoiding certain patterns. Before starting, we give a small bunch of
definitions that will be useful.

Let $Q,R$ be subsets of the ground set of a poset $\mathscr{P}$. We say that $Q$ is \emph{weakly below} $R$,
and write $Q\ngeq R$, whenever, for all $x\in Q$ and $y\in R$, we have $x\ngeq y$ in $\mathscr{P}$.
Given subsets $P_1 ,P_2 ,\ldots ,P_r$ of $\mathscr{P}$, we say that $(P_1 ,P_2 ,\ldots ,P_r)$ is a
\emph{weakly ordered partition} of $\mathscr{P}$ when it is a set partition of the ground set of $\mathscr{P}$
such that $P_1 \ngeq P_2 \ngeq \cdots \ngeq P_r$. Replacing every symbol $\ngeq$ with $\leq$
in the two above definitions, we obtain the definitions of ``$Q$ is \emph{below} $R$'' and
``$(P_1 ,P_2 ,\ldots ,P_r)$ is an \emph{ordered partition} of $\mathscr{P}$''.

Suppose that $\mathscr{P}$ and $\mathscr{Q}$ are two posets. The \emph{disjoint union}
$\mathscr{P}\cupdot \mathscr{Q}$ is the poset whose ground set is the disjoint union of the ground sets
of $\mathscr{P}$ and $\mathscr{Q}$ and such that $x\leq y$ in $\mathscr{P}\cupdot \mathscr{Q}$ whenever
either $x\leq y$ in $\mathscr{P}$ or $x\leq y$ in $\mathscr{Q}$. The \emph{linear sum}
$\mathscr{P}\oplus \mathscr{Q}$ is the poset whose ground set is the disjoint union of the ground sets
of $\mathscr{P}$ and $\mathscr{Q}$ and such that $x\leq y$ in $\mathscr{P}\oplus \mathscr{Q}$ whenever
either $x\leq y$ in $\mathscr{P}$ or $x\leq y$ in $\mathscr{Q}$ or else $x\in \mathscr{P}$ and
$y\in \mathscr{Q}$.

\begin{prop}
If $\mathscr{P}$ is the discrete poset of size $n$, then $SAv(\mathscr{P})$ is the class of all posets
of size $\leq n-1$.
\end{prop}

\emph{Proof.}\quad Clearly, any poset $\mathscr{Q}$ of size $\leq n-1$ strongly avoids $\mathscr{P}$,
since there cannot exist any injective function from $\mathscr{P}$ to $\mathscr{Q}$. Moreover,
if $\mathscr{Q}$ has size $\geq n$, then any injective map from $\mathscr{P}$ to $\mathscr{Q}$
trivially preserves the order ($\mathscr{P}$ does not have any order relation
among its elements, except of course the trivial ones coming from reflexivity), and so
$\mathscr{Q}\notin SAv (\mathscr{P})$.\cvd

\begin{prop}
If $\mathscr{P}$ is the poset of size $n$ containing a single covering relation, then
$SAv(\mathscr{P})$ contains all posets of size $\leq n-1$ and all (finite) discrete posets.
\end{prop}

\emph{Proof.}\quad This is essentially a consequence of the previous proposition. Indeed, every poset
$\mathscr{Q}$ of size $<n$ trivially strongly avoids $\mathscr{P}$; moreover, if $\mathscr{Q}$
has size $\geq n$ and strongly avoids $\mathscr{P}$, then it cannot contain any covering relation,
i.e. it is a discrete poset.\cvd

Recall that the \emph{height} of a poset is the maximum cardinality of a chain.

\begin{prop}\label{chain}
If $\mathscr{P}$ is the chain of size $n$, then $SAv(\mathscr{P})$ is the class of all finite posets
of height $\leq n-1$.
\end{prop}

\emph{Proof.}\quad If there is an injective function from $\mathscr{P}$ to a certain poset $\mathscr{Q}$
which preserves the order, then every pair of elements in the image of $\mathscr{P}$ must comparable,
i.e. $f(\mathscr{P})$ has to be a chain. This immediately yields the thesis.\cvd

Notice that, in this last case, that is when $\mathscr{P}$ is a chain, clearly
$SAv(\mathscr{P})=Av(\mathscr{P})$, since the up-set generated by $\mathscr{P}$ in $\mathfrak{X}_n$
consists of $\mathscr{P}$ alone.

A finite poset is called a \emph{flat} when it consists of a (possibly empty) antichain with an added maximum.

\begin{prop}
If $\mathscr{P}=$ \begin{tikzpicture}[scale=0.25] \draw (0,1) [fill] circle (.1);
\draw (1,0) [fill] circle (.1); \draw (2,1) [fill] circle (.1);
\draw (0,1) -- (1,0); \draw (1,0) -- (2,1); \end{tikzpicture} , then $SAv(\mathscr{P})$ is the class of
all disjoint unions of flats. As a consequence, $|SAv_n (\mathscr{P})|=p_n$, the number of integer partitions of $n$.
\end{prop}

\emph{Proof.}\quad Thanks to Proposition \ref{upset} and \ref{chain}, we observe that
$SAv(\mathscr{P})=Av($ \begin{tikzpicture}[scale=0.2] \draw (0,1) [fill] circle (.1);
\draw (1,0) [fill] circle (.1); \draw (2,1) [fill] circle (.1);
\draw (0,1) -- (1,0); \draw (1,0) -- (2,1); \end{tikzpicture} ,\, \begin{tikzpicture}[scale=0.2]
\draw (0,0) [fill] circle (.1); \draw (0,1) [fill] circle (.1); \draw (0,2) [fill] circle (.1);
\draw (0,0) -- (0,1); \draw (0,1) -- (0,2); \end{tikzpicture} $)$, and so in particular,
if $\mathscr{Q}\in SAv(\mathscr{P})$, then $\mathscr{Q}$ has height at most 1.
Moreover, any subset of cardinality 3 of $\mathscr{Q}$ cannot have minimum; thus, any three elements
in the same connected component are either an antichain or one of them is greater than the remaining two.
This means that each connected component of $\mathscr{Q}$ is a flat.

Concerning enumeration, the class of posets of size $n$ whose connected components are flats is in bijection
with the class of integer partitions of $n$: just map each of such posets into the integer partition
of the cardinality of its ground set whose parts are the cardinalities of the connected components (and observe that the order structure of a flat
is completely determined by its cardinality). From this observation the thesis follows.\cvd

We close this section with some general results which allow to understand the class $SAv(\mathscr{P})$
when $\mathscr{P}$ is built from simpler posets using classical operations.

\begin{prop}
Let $\mathscr{P},\mathscr{Q}$ be two posets.

\begin{enumerate}

\item $\mathscr{R}\in SAv(\mathscr{P}\cupdot \mathscr{Q})$ if and only if, for every partition $(R_1 ,R_2 )$
into two blocks of the ground set of $\mathscr{R}$, denoting with $\mathscr{R}_1 ,\mathscr{R}_2$
the associated induced subposets, $\mathscr{R}_1 \in SAv(\mathscr{P})$ or
$\mathscr{R}_2 \in SAv(\mathscr{Q})$.

\item If $\mathscr{R}\in SAv(\mathscr{P}\oplus \mathscr{Q})$, then, for every ordered partition
$(R_1 ,R_2 )$ into two blocks of $\mathscr{R}$, $\mathscr{R}_1 \in SAv(\mathscr{P})$ or
$\mathscr{R}_2 \in SAv(\mathscr{Q})$. If $\mathscr{R}\notin SAv(\mathscr{P}\oplus \mathscr{Q})$,
then there exists a weakly ordered partition $(R_1 ,R_2 )$ into two blocks of $\mathscr{R}$
such that $\mathscr{R}_1 \notin SAv(\mathscr{P})$ and $\mathscr{R}_2 \notin SAv(\mathscr{Q})$.

\end{enumerate}

\end{prop}

\emph{Proof.}

\begin{enumerate}

\item An occurrence of $\mathscr{P}\cupdot \mathscr{Q}$ in $\mathscr{R}$ consists of an occurrence of $\mathscr{P}$
and an occurrence of $\mathscr{Q}$ whose ground sets are disjoint and with no requirements about the order relations
among pairs of elements $(x,y)$ such that $x\in \mathscr{P}$ and $y\in \mathscr{Q}$. Therefore,
if $\mathscr{R}\in SAv(\mathscr{P}\cupdot \mathscr{Q})$ and $(R_1 ,R_2 )$ is a partition of the ground
set of $\mathscr{R}$, then it is clear that, if $\mathscr{R}_1$ weakly contains $\mathscr{P}$, then necessarily
$\mathscr{R}_2$ strongly avoids $\mathscr{Q}$. Vice versa, if $\mathscr{R}$ weakly contains $\mathscr{P}\cupdot \mathscr{Q}$,
then clearly there exists an occurrence of $\mathscr{P}$ whose complement weakly contains $\mathscr{Q}$.

\item An occurrence of $\mathscr{P}\oplus \mathscr{Q}$ in $\mathscr{R}$ consists of an occurrence of $\mathscr{P}$
and an occurrence of $\mathscr{Q}$ whose ground sets are disjoint and such that every element of $\mathscr{P}$
is less than every element of $\mathscr{Q}$. Thus, if $\mathscr{R}\in SAv(\mathscr{P}\oplus \mathscr{Q})$ and $(R_1 ,R_2 )$
is an ordered partition of $\mathscr{R}$ such that $R_1$ weakly contains $\mathscr{P}$, then necessarily $R_2$
strongly avoids $\mathscr{Q}$, since the ground set of $\mathscr{P}$ lies below $R_2$. On the other hand,
if $\mathscr{R}$ weakly contains $\mathscr{P}\oplus \mathscr{Q}$, then the partition $(R_1 ,R_2  )$ of $\mathscr{R}$
in which $R_1$ is the down-set generated by an occurrence of $\mathscr{P}$ (and, of course, $R_2$ is the
complement of $R_1$, and so an up-set) is a weakly ordered partition having the required properties.\cvd

\end{enumerate}

Another simple, general result involving the disjoint union of posets is the following.

\begin{prop} If $\mathscr{Q}\in SAv(\mathscr{P})$ and $\mathscr{P}$ is connected,
then $\mathscr{Q}$ is the disjoint union of a family of posets strongly avoiding $\mathscr{P}$.
\end{prop}

\emph{Proof.}\quad Indeed, take $\mathscr{Q}\in SAv(\mathscr{P})$ and suppose that $\mathscr{Q}$ is not connected
(otherwise the thesis is trivial). Since $\mathscr{P}$ is connected,
any occurrence of $\mathscr{P}$ in $\mathscr{Q}$ would be connected too (since such an occurrence is essentially
$\mathscr{P}$ with possibly some added order relations), so, in order $\mathscr{Q}$ to strongly avoid $\mathscr{P}$,
each connected component of $\mathscr{Q}$ has to strongly avoid $\mathscr{P}$, as desired.\cvd

\subsection{How Schr\"oder tableaux come into play}

Our introduction of weak poset patterns is motivated by the role they have in the description of Schr\"oder tableaux.
Let $S$ be a Schr\"oder tableau, and let $\mathcal{I}_S$ be the associated set of intervals, as defined in Subsection
\ref{interval}. Consider the interval order associated with $\mathcal{I}_S$, to be denoted $\mathcal{I}_S$ as well.
The map $S\mapsto \mathcal{I}_S$ from Schr\"oder tableaux to interval orders is clearly neither injective nor surjective.
Therefore two natural questions concerning such a map arise.

\begin{enumerate}

\item Given an interval order $\mathcal{I}$, does there exist a Schr\"oder tableau $S$ such that
$\mathcal{I}=\mathcal{I}_S$?

\item In case of a positive answer to the previous question, how many Schr\"oder tableaux associated
with a given interval order are there?

\end{enumerate}

Below we give an answer to the first question. Recall that $\mathbf{N}$ denotes the set of natural numbers (without $0$),
which will be endowed with its usual total order structure.

\begin{teor}
Let $\mathcal{I}$ be an interval order of size $n$. There exists a Schr\"oder tableau $S$ such that $\mathcal{I}=\mathcal{I}_S$
if and only if $\mathcal{I}$ weakly contains a down-set of size $n$ of $\mathbf{N}\times \mathbf{N}$.
\end{teor}

\emph{Proof.}\quad Suppose first that $\mathcal{I}=\mathcal{I}_S$, for some Schr\"oder tableau $S$. This means that
we can use the map $S\mapsto \mathcal{I}_S$ to label the set $C$ of all pairs of twin cells and of all lonely cells of $S$
with the elements of $\mathcal{I}$, so that $C$ can be identified with $\mathcal{I}$.
Consider the function $f:\mathcal{I}\rightarrow \mathbf{N}\times \mathbf{N}$ mapping $I$ to the pair
$(n_I ,m_I )$ such that $n_I$ (resp. $m_I$) is the row (resp. column)
of the pair of twin cells or of the lonely cell associated with $I$ (as usual, rows and columns are enumerated
from top to bottom and from left to right, respectively). We show that $f(\mathcal{I})$ is a down-set of $\mathbf{N}\times \mathbf{N}$:
indeed, if $I\in \mathcal{I}$ and $(n,m)\in \mathbf{N}\times \mathbf{N}$ are such that $(n,m)\leq f(I)=(n_I ,m_I )$,
then the tableau $S$ has at least $n_I$ rows and $m_I$ columns, so in particular it exists a pair of twin cells or a lonely cell
at the crossing of row $n$ and column $m$ (since $n\leq n_I$ and $m\leq m_I$); denoting with $J$ the associated interval of
$\mathcal{I}$, we then have that $(n,m)=f(J)\in f(\mathcal{I})$, as desired. By construction $f$ is injective,
since two distinct intervals of $\mathcal{I}$ correspond to two distinct cells of $S$, which of course cannot lie
both in the same row and in the same column of $S$. We can thus consider the inverse $g:f(\mathcal{I})\rightarrow \mathcal{I}$
of $f$ on $f(\mathcal{I})$. We now show that $g$ is order-preserving: indeed, consider $I,J\in \mathcal{I}$ and suppose that
$(n_I ,m_I )=f(I)\leq f(J)=(n_J ,m_J )$; the tableau $S$ has a pair of twin cells or a lonely cell at the crossing of
row $n_J$ and column $m_I$ (since $m_I \leq m_J$), and the associated interval, call it $K$, is such that $K\leq J$ in
$\mathcal{I}$; moreover, since $n_I \leq n_J$, it is not difficult to realize that $I\leq K$, hence we conclude that $I\leq J$.
Therefore we have proved that $f(\mathcal{I})$ is a down-set (of size $n$) of $\mathbf{N}\times \mathbf{N}$ and that
$g:f(\mathcal{I})\rightarrow \mathcal{I}$ is an injective order-preserving map, i.e. $f(\mathcal{I})$ is a weak pattern of
$\mathcal{I}$: this is exactly the thesis.

In the other direction, suppose that $\mathcal{I}$ weakly contains a down-set $D$ of $\mathbf{N}\times \mathbf{N}$ of size $n$.
In other words, $D$ is a coarsening of $\mathcal{I}$ (i.e. it is obtained from $\mathcal{I}$ by possibly removing
some order relations only, leaving untouched its ground set). This means that there is an order-preserving injective map
$g:D\rightarrow \mathcal{I}$. As we have already recalled, it is not restrictive to suppose that
the endpoints of the elements of $\mathcal{I}$ are all distinct and constitute an initial segment of $\mathbf{N}$.
Consider the tableau $S$ having a pair of twin cells in row $n$ and column $m$ if and only if $(n,m)\in D$ and, in such case,
the two cells are filled in with the endpoints of the interval $g(n,m)\in \mathcal{I}$. It is clear that $S$ is a Schr\"oder tableau:
since $D$ is a down-set, rows are left-justified; moreover rows and columns are increasing because $g$ is order-preserving.
Now it is not difficult to realize that $\mathcal{I}=\mathcal{I}_S$, which is precisely what we wanted to show.\cvd

\emph{Remark.}\quad Notice that, in the second part of the above proof, we construct a Schr\"oder tableau $S$
\emph{having no lonely cells} such that $\mathcal{I}=\mathcal{I}_S$ . This again shows the fact that
there can be several Schr\"oder tableaux associated with the same interval order.

\section{Further work}\label{conclusions}

The algebraic and combinatorial properties of the distributive
lattice \Sch of Schr\"oder shapes needs to be further
investigated. In particular, the analogies with differential
posets should be much deepened, for instance trying to understand
the role of the lowering and raising operators (a fundamental tool for
computations in differential posets) in the Schr\"oder lattice, or even
in the more general setting of $\varphi$-differential posets.

\bigskip

%We have shown that Schr\"oder partitions are the fixed points of the map $c_2 ^2$.
%Can we characterize in an analogous way the set of fixed points of $c_n ^2$,
%for $n\geq 3$?
%
%\bigskip

We have just started the characterization and enumeration of
permutations having a given Schr\"oder insertion tableau. Many
more shapes should be investigated. Moreover, we still have to
understand the role of the recording tableau.

\bigskip

Can we find a nice closed formula for the number of Schr\"oder tableaux
of a given shape? In the case of Young tableaux there is a famous
\emph{hook formula}, which however seems to be unlikely in our
case, since we have numerical evidence that, for certain shapes,
this number has large prime factors.

\bigskip

The alternative presentation of Schr\"oder tableaux in terms of interval orders
and weak poset patterns might have more secrets to reveal. For instance,
the enumeration of Schr\"oder tableaux associated with a given interval order
is entirely to be done. In a different direction, the topic of strong pattern avoidance
for posets seems to be an interesting line of research in its own,
independently from its relationship with Schr\"oder tableaux.

\bigskip

The analogies between Young tableaux and Schr\"oder tableaux
should be investigated more, especially from a purely
algebraic point of view. Combinatorial objects related to Young
tableaux, such as Schur functions and the plactic monoid, as well
as algorithmic and algebraic constructions, such as
Sch\"utzenberger's \emph{jeu de taquin} \cite{Sa}, the Littlewood-Richardson
rule and the Schubert calculus on Grassmannians and flag
varieties, could have some interesting counterparts in the context
of Schr\"oder tableaux.

%\bigskip
%
%\textbf{Acknowledgement.}\quad We are very grateful to the anonymous referees for an extremely careful reading and many
%useful suggestions. This has led to the correction of some mistakes and to a significant improvement of the presentation.


\begin{thebibliography}{99}

\bibitem[AA]{AA}
Albert, M., Atkinson, M. D.:
Pattern classes and priority queues.
Pure Math. Appl. (PU.M.A.) 23, 161--177 (2012).

\bibitem[A]{A}
Andrews, G. E.:
Partition identities.
Adv. Math. 9, 10--51 (1972).

\bibitem[Be]{Be}
Bergeron, F.:
Algebraic combinatorics and coinvariant spaces.
CMS Treatise in Mathematics, A K Peters/CRC Press, 2009.

\bibitem[Bo]{Bo}
B\'ona, M.:
Combinatorics of permutations.
Discrete Mathematics and its Applications, CRC Press, Taylor \& Francis Group,
2012.

\bibitem[BCDK]{BCDK}
Bousquet-Melou, M., Claesson, A., Dukes, M., Kitaev, S.:
$(\mathbf{2}+\mathbf{2})$-free posets, ascent sequences and pattern avoiding permutations.
J. Combin. Theory Ser. A 117, 884-909 (2010).

\bibitem[D]{D}
Drake, B.:
Limits of areas under lattice paths.
Discrete Math. 309, 3936--3953 (2009).

\bibitem[FP]{FP}
Ferrari, L., Pinzani, R.:
Lattices of lattice paths.
J. Statist. Plann. Inference 135, 77--92 (2005).

\bibitem[Fe]{Fe}
Ferrari, L.:
Schr\"oder partitions and Schr\"oder tableaux.
Lecture Notes in Computer Science 9538,
Combinatorial Algorithms (Zsuzsanna Liptak, William F. Smyth, eds.), 161--172 (2016).

\bibitem[Fi]{Fi}
Fishburn, P. C.:
Intransitive indifference with unequal indifference intervals.
J. Math. Psych. 7, 144--149 (1970).

\bibitem[FRT]{FRT}
Frame, J. S., Robinson, G. de B., Thrall, R. M.:
The hook graphs of the symmetric group.
Canad. J. Math. 6, 316--325 (1954).

\bibitem[GL]{GL}
Griggs, J. R, Lu, L.:
On families of subsets with a forbidden subposet.
Combin. Probab. Comput. 18, 731--748 (2009).

\bibitem[HM]{HM}
Heubach, S., Mansour, T.:
Combinatorics of Compositions and Words.
Chapman \& Hall/CRC, Taylor \& Francis Group, Boca Raton, London, New York, 2009.

\bibitem[KT]{KT}
Katona, G. O., Tarjan, T. G.:
Extremal problems with excluded subgraphs in the $n$-cube.
In: Graph Theory, pp. 84--93, Springer, 1983.

\bibitem[Ki]{Ki}
Kitaev, S.:
Patterns in permutations and words.
EATCS Monographs in Theoretical Computer Science, Springer-Verlag, 2011.

\bibitem[Kn]{Kn}
Knuth, D. E.:
Permutations, matrices, and generalized Young tableaux.
Pacific J. Math. 34, 709--727 (1970).

\bibitem[LR]{LR}
Littlewood, D. E., Richardson, A. R.:
Group characters and algebra.
Philos. Trans. R. Soc. Lond. Ser. A Math. Phys. Eng. Sci. 233, 99--141 (1934).

\bibitem[M]{M}
Mansour, T.:
Combinatorics of Set Partitions.
Chapman \& Hall/CRC, Taylor \& Francis Group, Boca Raton, London, New York, 2012.

\bibitem[R]{R}
Robinson, G. de B.:
On the representations of the symmetric group.
Amer. J. Math. 60, 745--760 (1938).

\bibitem[Sa]{Sa}
Sagan, B. E.:
The Symmetric Group: Representations, Combinatorial Algorithms, and Symmetric Functions.
Graduate Texts in Mathematics 203 (2nd ed.), New York: Springer, 2001.

\bibitem[Sc]{Sc}
Schensted, C.:
Longest increasing and decreasing subsequences.
Canad. J. Math. 13, 179--191 (1961).

\bibitem[SiSc]{SiSc}
Simion, R., Schmidt, F. W.:
Restricted permutations.
European J. Combin. 6, 383--406 (1985).

\bibitem[Sl]{Sl}
Sloane, N. J. A.:
The On-Line Encyclopedia of Integer Sequences.
http://oeis.org.

\bibitem[St]{St}
Stanley, R. P.:
Differential posets.
J. Amer. Math. Soc. 1, 919--961 (1988).

\end{thebibliography}
\end{document}